\documentclass[a4paper]{article}

\usepackage{amsmath}

\usepackage{amsthm}
\usepackage[dvips]{graphicx}

\usepackage{epsfig}

\usepackage{amssymb}
\usepackage{latexsym}
\newtheorem{Definition}{Definition}
\newtheorem{Theorem}{Theorem}
\newtheorem{Proposition}{Proposition}
\newtheorem{Lemma}{Lemma}
\newtheorem{Corollary}{Corollary}
\newtheorem{Remark}{Remark}

\def\L{{\mathcal L}}
\def\K{{\mathcal K}}
\def\O{{\mathcal O}}

\def\ZZ{{\mathbf Z}}
\def\QQ{{\mathbf Q}}

\def\sph{{\bf S}^3}

\begin{document}

\title{Hyperbolic isometries versus symmetries of links}

\author{Luisa Paoluzzi and Joan Porti\footnote{Partially supported by the 
Spanish MEC through grant BFM2003-03458 and European funds FEDER}} 
\date{\today}

\maketitle

\begin{abstract}
We prove that every finite group is the orientation-preserving isometry group 
of the complement of a hyperbolic link in the $3$-sphere.
\vskip 2mm


\vskip 2mm


\end{abstract}

\section{Introduction}
\label{section:intorduction}

There are several well-known results in the study of finite group actions on 
special $3$-manifolds, starting with the work \cite{Mi} of Milnor, who gave a 
list of all finite groups which are susceptible to act freely and 
orientation-preserving on spheres. In the $3$-dimensional case, by
elliptization of three-manifolds (see \cite{P1}, \cite{P2}, \cite{P3}, 
\cite{CM} and \cite{B}) and Thurston's orbifold theorem (see \cite{T} and
\cite{BLP}, \cite{BP} and \cite{CHK}), the only finite groups which 
can act on $\sph$ preserving the orientation are the finite subgroups of 
$SO(4)$. It happens also that the list of finite groups acting on integral and, 
more generally, $\ZZ/2$-homology spheres is quite restricted \cite{MZ}.

%
%

\begin{Definition}
A \emph{group of symmetries} of a (non-oriented) link $L$ in $\sph$ is a finite 
group $G$ acting on $\sph$ preserving the orientation and leaving $L$ 
invariant.
\end{Definition}

In particular, a group of symmetries is a finite subgroup of $SO(4)$.
%
%
%
Indeed, if $L$ is a non-trivial knot, it follows from Smith's conjecture 
\cite{MB} that the only possible groups of symmetries for $L$ are either cyclic
or dihedral. Note also that, in general, a link does not have a unique 
(maximal) group of symmetries \cite{S} (up to conjugation), but this is indeed 
the case if the link is hyperbolic. Observe in fact, that each symmetry of the 
link $L$ induces an orientation-preserving diffeomorphism of the complement of 
the link, which on its turn, if $L$ is hyperbolic, gives rise to an isometry of 
the hyperbolic structure. (For basic facts and definitions in hyperbolic
geometry the reader is referred to \cite{R}.)

One can then ask: \emph{What is the relation between symmetries of a hyperbolic 
link and isometries of its complement?} As knots are determined by their 
complements \cite{GL},  isometries of the complement of a knot are also 
symmetries. However, when the link has several components, the group of 
symmetries of a hyperbolic link is only a subgroup of the group of isometries 
of its complement. In general this subgroup can be proper, for a generic 
isometry does not need  to preserve a peripheral structure on the cusps. 

Examples of isometries of exteriors of hyperbolic links which are not induced
by symmetries can be found in \cite{HW}, where both groups are computed 
using Jeff Week's program SnapPea for links with at most nine crossings. For 
most links in their list both groups coincide, and when they do not, the 
index of the symmetry group in the isometry one is rather small ($\le4$). Also, 
the symmetry groups which appear are of very special types (either abelian, or 
dihedral). 

This behavior is certainly due to the limited number of crossings and
components considered. Indeed, in this paper we prove:

\medskip

\begin{Theorem}\label{thm}
Every finite group $G$ is the group of orientation preserving isometries of the
complement of some hyperbolic link in $\sph$. Moreover, $G$ acts freely.
\end{Theorem}


It was already known that every finite group can be realized as an isometry 
group of some closed hyperbolic $3$-manifold: Kojima proved that every finite 
group can be realized as the full group of isometries of a closed hyperbolic 
manifold \cite{K}, with a not necessarily free action. Besides,  Cooper and 
Long showed that every finite group acts freely on some hyperbolic rational 
homology sphere \cite{CL}; in this case the full isometry group might a priori
be larger.

Although the aforementioned results are related to the one presented in this
paper, it deserves to be stressed that there are some peculiarities which come
from the fact that we want to control several things at the same time: the
structure of the manifold (i.e.\ the complement of a link in $\sph$), the type 
of action (i.e.\ free), and the fact that no extra isometries are introduced. 

Given a finite group $G$, it is not difficult to construct a link in $\sph$ 
(and even a hyperbolic one, thanks to a result of Myers \cite{My1}) whose
complement admits an effective free $G$-action (see Proposition~\ref{prop}), 
but it is rather delicate to ensure that $G$ coincides with the whole isometry
group of the complement. Notice that the naive idea to drill out some
``asymmetric" link in a $G$-equivariant way does not work, for the fact of 
removing simple closed curves has the effect of possibly increasing the group 
of isometries.

The idea is thus to choose the link in such a way that its complement contains
some very rigid structure (in our case, totally geodesic pants) and use it
to control the isometry group. 

As a by-product of our result we are able to prove that Cooper and Long's
result can be rigidified, that is, that the hyperbolic rational homology sphere
on which the group $G$ acts freely can be chosen so that $G$ coincides with the 
full orientation preserving isometry group.

\begin{Corollary}\label{cor}
Any finite group $G$ is the full orientation-preserving isometry group of a
hyperbolic rational homology sphere. Moreover $G$ acts freely.

\end{Corollary}

Before passing to the proof of the result, which will be the content of
Sections~\ref{section:construction} and \ref{section:nomoreisos},
we state the following:

\begin{Proposition}\label{linksymmetries}
For any finite group $H$ acting on $\sph$, there exists a hyperbolic link
whose group of symmetries is precisely $H$.
\end{Proposition}

\begin{proof}
We start by noticing that Myer's theorem works in an $H$-invariant setting 
(it suffices to choose an $H$-invariant triangulation) if one does not require 
the resulting hyperbolic link to be connected (see the proof of
Lemma~\ref{claim:asymmetric hyperbolic link} for a similar reasoning). In this
case one can indeed exploit the naive idea to add some ``asymmetric"
but $H$-invariant link to ensure that $H$ coincides with the whole group of
symmetries. More precisely, given any finite group $G$ acting on $\sph$, take 
$|H|$ copies of a knot $K$ which admits no symmetry, in such a way that each 
copy is contained in a ball which misses all the fixed-point sets of elements
of $H$, and such that the pairs (ball,$K$) are freely permuted by the action of 
$H$. It is now possible, using a result of Myers \cite{My2}, to find an 
$H$-invariant link $\Lambda$ so that no component of $\Lambda$ is the knot $K$ 
and the complement of $L=\Lambda\cup G\cdot K$ is hyperbolic. It is now 
straightforward to see that the group of symmetries of the link $L$ is 
precisely $H$.
\end{proof}

Note finally that in both constructions, Theorem~\ref{thm} and 
Proposition~\ref{linksymmetries}, the number of components of the links 
constructed is pretty large. It would be interesting to understand whether, for 
a given group $G$, there is some bound on the number of components of a link on 
which $G$ acts as group of symmetries/isometries and under which conditions the 
lower bound can be realized.


\section{Construction of the link}
\label{section:construction}

The following result is easy and is the starting point to prove our main 
result.

\medskip

\begin{Proposition}\label{prop}
Every finite group acts effectively and freely by orientation-preserving
isometries on the complement of some hyperbolic link in the $3$-sphere.
\end{Proposition}

\begin{proof}
Fix a finite group $G$. It is not difficult to find a closed $3$-manifold $M$
on which $G$ acts effectively and freely; moreover, Cooper and Long \cite{CL} 
showed that $M$ can be chosen to be hyperbolic and a rational homology sphere.  
The Dehn surgery theorem of Lickorish \cite{L} and Wallace \cite{W} assures 
the existence of a link $\L$ contained in $M$ whose exterior is contained in 
the $3$-sphere, i.e. $M\setminus\L=\sph\setminus L$, where $L$ is a link. We 
want to show that $\L$ can be chosen $G$-invariant, thanks to a general 
position argument. Take the quotient $(M,\L)/G$. By perturbing slightly the 
image of $\L$ inside the quotient, we can assume that it has no 
self-intersection. Note that the complement of the $G$-invariant link of $M$ 
obtained this way is again contained in $\sph$, for we can assume that the 
perturbation performed did not affect the isotopy class of the components 
corresponding to the original link. Now, using a result of Myers \cite{My1}, we 
can find in $(M,\L)/G$ a knot $\K$ whose exterior is hyperbolic. The preimage 
in $M$ of $\L/G\cup\K$ is a link with the desired properties.
\end{proof}

\bigskip

\begin{figure}[ht]
\begin{center}
\includegraphics[width=6cm]{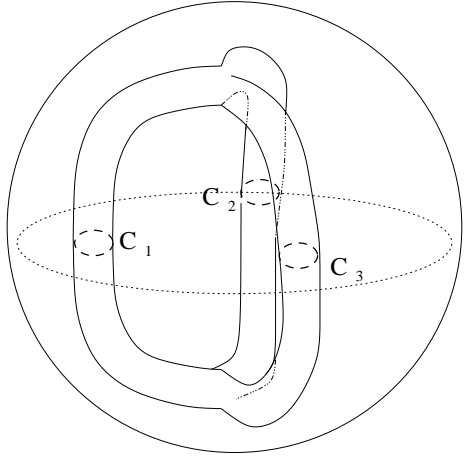}
\end{center}
\caption{The handlebody $H$ contained in a ball and the three curves $C_i$.}
\label{handle}
\end{figure}

\noindent\emph{Proof of Theorem~\ref{thm}.}
We start with $G$ acting on $\sph\setminus L$ for some link $L\subset \sph$, 
as in Proposition~\ref{prop}. A priori the group 
$\operatorname{Isom}^+(\sph\setminus L)$ can be larger than $G$, and we shall 
modify the link so that both groups are precisely the same. Denote
by $N$ the quotient $(\sph\setminus L)/G$. Choose a genus-$2$ unknotted 
handlebody contained in a ball of $N$. Using \cite{My1}, we can find a 
hyperbolic knot $\K$ in $N\setminus\mathring{H}$, whose exterior is, moreover, 
anannular. Fix three meridional curves $C_i$, $i=1,2,3$, on the boundary of $H$ 
as shown in Figure~\ref{handle}. Observe that each $C_i$ is non-separating, but 
the three of them cut $\partial H$ in two pairs of pants. The manifold
$(N\setminus\mathring{H})\setminus(\K\cup\bigcup_{i=1}^3 C_i)$ has a 
non-compact boundary, consisting of two cusped pairs of pants, and admits a 
hyperbolic structure with totally geodesic boundary: this last fact follows 
from Thurston's hyperbolization theorem, since the manifold is atoroidal and 
its compact core is anannular. Notice that the curves $C_i$ correspond to rank 
one cusps.

We consider the handlebody $H$.

\begin{Lemma}\label{claim:asymmetric hyperbolic link}
There exists a link $\Lambda$ inside $H$
such that
$H\setminus(\Lambda\cup\bigcup_{i=1}^3 C_i)$ satisfies the following 
properties:
\begin{enumerate}

\item It is hyperbolic with totally geodesic boundary.

\item It contains a unique geodesic $\eta$ of minimal length.

\item Its orientation preserving isometry group is trivial.

\end{enumerate}
\end{Lemma}

\begin{proof}
Consider the group $\ZZ/2\times{\bf D}_3$, where ${\bf D}_3$ denotes the 
dihedral group of order $6$, which coincides with the symmetric group on three 
elements. This group acts on the pair $(H,\bigcup_{i=1}^3 C_i)$ as illustrated 
in Figure~\ref{sym}a: the central element of order two is the hyperelliptic 
involution fixing all three curves $C_i$ while exchanging the two pairs of 
pants (its ``axis" is a circle in the figure), the elements of ${\bf D}_3$ 
leave invariant each pair of pants and permute the curves $C_i$. Consider now 
the quotient orbifold (see Figure~\ref{sym}b)
$$
\O=(H\setminus\bigcup_{i=1}^3 C_i) \Big/  \ZZ/2\times{\bf D}_3.
$$
Remove from $\O$ a simple closed curve 
$\gamma$ contained in a ball which does not meet the singular locus. Choose a 
triangulation of the compact core of $\O\setminus\gamma$ which contains the 
singular locus: such triangulation lifts to a 
$\ZZ/2\times{\bf D}_3$-equivariant triangulation of the compact core of 
$H\setminus(\tilde\gamma\cup\bigcup_{i=1}^3 C_i)$, where 
$\tilde\gamma$ denotes the lift of $\gamma$. Taking the second baricentric 
subdivision of such triangulation, one can get a ``special handle 
decomposition" and the very same proof as in \cite[Theorem 6.1]{My1} shows that 
one can find a hyperbolic link
$
\Lambda_0\subset H\setminus(\tilde\gamma\cup\bigcup_{i=1}^3 C_i)
$
 which is $\ZZ/2\times{\bf D}_3$-invariant and whose exterior 
$$
W=  H\setminus(\tilde\gamma\cup\Lambda_0\cup \bigcup_{i=1}^3 C_i)
$$
 has 
anannular compact core. 
Consider now $\tilde\gamma$: it 
consists of twelve connected components, for each of which we choose a 
meridian-longitude system $(\mu,\lambda)$. One can perform hyperbolic Dehn 
surgery with meridian curves $\mu+n\,\lambda$, $n\gg 1$ so that the lengths of 
the surgered geodesics are pairwise distinct and of shortest length inside the 
resulting manifold $H\setminus(\Lambda\cup\bigcup_{i=1}^3 C_i)$, $\Lambda$ 
being the image of $\Lambda_0$ after surgery.

\begin{figure}
\begin{center}
\includegraphics[width=10cm]{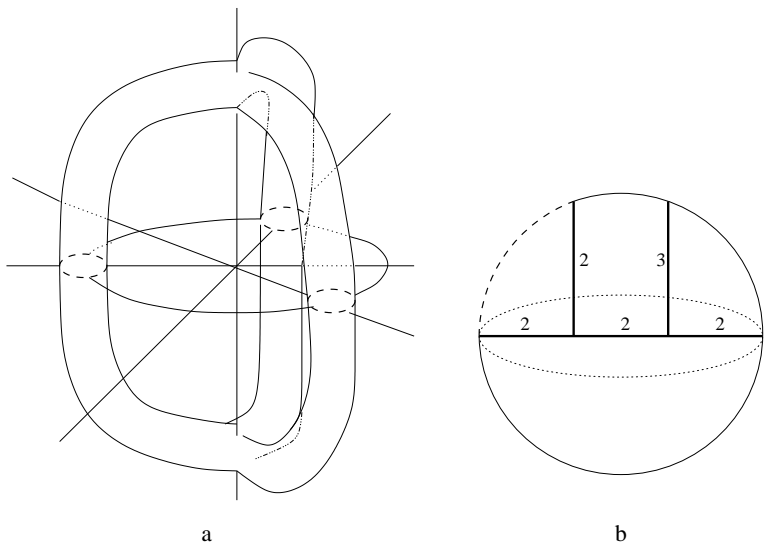}
\end{center}
\caption{The group $\ZZ/2\times{\bf D}_3$ acting on $H$ and the quotient
orbifold $\O$.}
\label{sym}
\end{figure}

Consider now an orientation-preserving isometry of 
$H\setminus(\Lambda\cup\bigcup_{i=1}^3 C_i)$. Since they have minimal length, 
the twelve geodesics obtained by hyperbolic surgery must be left invariant by 
the isometry, which thus induces an isometry of the exterior of the geodesics. 
This means that the isometry must act as one of the elements of 
$\ZZ/2\times{\bf D}_3$, since it is determined by its action on the boundary 
and since $\ZZ/2\times{\bf D}_3$ is the complete group of positive isometries 
of the boundary. But only the identity element of $\ZZ/2\times{\bf D}_3$ 
extends to $H\setminus(\Lambda\cup\bigcup_{i=1}^3 C_i)$, for all twelve 
geodesics must be left setwise fixed, because their lengths are pairwise 
distinct. 
\end{proof}

The boundary components of $H\setminus(\Lambda\cup\bigcup_{i=1}^3 C_i)$
are two totally geodesic pairs of pants with cusps (topologically, they are 
$\partial H\setminus \bigcup_{i=1}^3 C_i$). By abuse of notation, $C_i$ denotes 
the curve and also the corresponding cusp.

\begin{Lemma}\label{Lemma:no_pants}
The only embedded geodesic pair of pants in 
$H\setminus(\Lambda\cup\bigcup_{i=1}^3 C_i)$ having cusps $C_1$, $C_2$ and 
$C_3$ are its boundary components.
\end{Lemma}

\begin{proof}  
Let $P_{}$ be an embedded pair of pants with cusps $C_1$, $C_2$ and $C_3$. By 
the geometry of Margulis tubes, $P_{}$ must be disjoint from the cores of the 
filling tori, thus $P_{}\subset W$ (recall that $W$ was defined in the proof of 
Lemma~\ref{claim:asymmetric hyperbolic link}, and that it is the result of 
removing the twelve shortest geodesics). The pants $P_{}$ can be rendered 
totally geodesic also in $W$, because $W$ minus an open tubular neighborhood of 
$P_{}$
$$
W\setminus\mathcal N (P_{})
$$
has still an irreducible, atoroidal, and  anannular compact core. (This follows 
from the fact that every simple closed curve in $P_{}$ is either compressible 
or boundary parallel.) By Thurston's hyperbolization,  $W\setminus \mathcal 
N(P_{})$ is hyperbolic with totally geodesic boundary, consisting of two copies 
of $P$ that we can glue back by an isometry.

Next we claim that $P_{}$ is equivariant by the action of 
$\ZZ/2\times{\bf D}_3$. To prove it, we look carefully at the geometry of the 
cusps $C_i$. A horospherical section of the cups $C_i$ is a Euclidean annulus 
$A_i$, and $P_{}$ intersects the annulus $A_i$ in a circle metrically parallel 
to the boundary components $\partial A_i$. Thus the intersection of $P_{}$ with 
those annuli must be equivariant. Since the action of $\ZZ/2\times{\bf D}_3$ on 
$P$ is determined by its restriction to the cusps, it follows that $P$ is 
equivariant.


\begin{figure}[ht]
\begin{center}
\includegraphics[width=6cm]{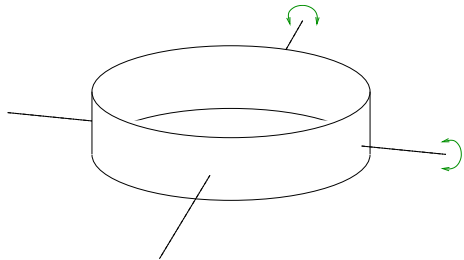}
\end{center}
\caption{The action of $\ZZ/2\oplus \ZZ/2$ on the annulus $A_i$.}
\label{annulus}
\end{figure}

Now we consider the possible quotients of $P_{}$ by its stabilizer, and look at 
how they project on $W/(\ZZ/2\times{\bf D}_3)$. We remark that, on the circle 
of the intersection of $P_{}$ with each horospherical annulus $A_i$, the 
projection to the quotient restricts to a map of the circle $P\cap A_i$ to its 
quotient which is either $2$  to $1$, or $4$ to $1$ (Figure~\ref{annulus}). In 
the former case, the quotient is a circle, in the latter, an interval with 
mirror boundary. This implies that the stabilizer of $P_{}$ must contain at 
least three orientation preserving involutions, one for each cusp, and 
corresponding to the product of the two involutions whose axes are pictured in
Figure~\ref{annulus} (i.e.\ a rotation on the $S^1$ factor).
 As the orientation preserving isometry group of the pants is $\mathbf D_3$, it follows that the only possible quotients for $P_{}$ are 
the ones described in Figure~\ref{quotients}.

\begin{figure}[ht]
\begin{center}
\includegraphics[width=6cm]{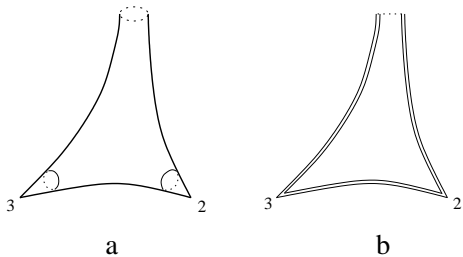}
\end{center}
\caption{The possible quotients of $P_{}$. Case $b$ is non orientable, has 
mirror points (double lines), and corners.}
\label{quotients}
\end{figure}

In fact case b in Figure~\ref{quotients} is not possible, as the mirror points 
have to agree with singular points of order $2$ in the orbifold $\mathcal O$ in 
Figure~\ref{sym}b. 

We claim that the quotient of $P$ is parallel to the boundary in $\O$. 
To see that, notice that $W/(\ZZ/2\times{\bf D}_3)$ is
the exterior in $\O$ 
 of an anannular hyperbolic link with two components, $\gamma$
and the quotient of $\Lambda_0$.
Those components cannot be separated by the quotient of $P$, because otherwise
the disc in $\O$ bounded by $\gamma$ would give either a compressing disc
or an essential annulus in $W\setminus \mathcal N(P)$, contradicting that $P$ 
is totally geodesic. Hence $\gamma$ and the quotient of $\Lambda_0$ are not 
separated by the quotient of $P$,  and therefore this quotient either
bounds a handlebody or is parallel to the boundary. It cannot bound a 
handlebody, since it is totally geodesic. Hence $P_{}$ must be a boundary 
component, because two parallel totally geodesic submanifolds must be the same.
\end{proof}

Consider now the manifold $N\setminus(\K\cup\Lambda\cup\bigcup_{i=1}^3
C_i)$. This is a hyperbolic manifold obtained by gluing together two hyperbolic
manifolds along their totally geodesic boundaries, which consist of two pairs 
of pants:
$$
N\setminus(\K\cup\Lambda\cup\bigcup_{i=1}^3 C_i)= 
\Big( H\setminus(\Lambda\cup\bigcup_{i=1}^3 C_i)\Big)
\cup
\Big( N\setminus(\mathring H\cup \K\cup\bigcup_{i=1}^3 C_i)  \Big).
$$ 
Since the hyperbolic structure of a pair of pants (with three cusps) is
unique, both pairs of pants are compatible when gluing and they remain 
totally geodesic inside the resulting manifold. Take now the lift of 
$N\setminus(\K\cup\Lambda\cup\bigcup_{i=1}^3 C_i)$ to the $3$-sphere: we obtain 
the complement of a hyperbolic link $L$ in $\sph$ on which $G$ acts freely, in 
particular $G$ is a subgroup in $\operatorname{Isom}^+(\sph\setminus L)$.

Notice that we glue along totally geodesic pants, therefore we do not make any 
deformation to the structure of each piece, and we get:

\begin{Remark}
The shortest geodesic $\eta$ in $H\setminus(\Lambda\cup\bigcup_{i=1}^3 C_i)$
from Lemma~\ref{claim:asymmetric hyperbolic link} may also be chosen to be the 
shortest in $N\setminus(\K\cup\Lambda\cup\bigcup_{i=1}^3 C_i)$, and its lift in 
$\sph\setminus L$ minimizes also the length spectrum.
\end{Remark}

\medskip

\section{Showing that there are no more isometries}
\label{section:nomoreisos}
 
It remains to show that $G=\operatorname{Isom}^+(\sph\setminus 
L)$. Arguing by contradiction, let $\varphi$ be an isometry of the complement 
of $L$ which is not in $G$. Consider now $\tilde\eta$ the lift in 
$\sph\setminus L$ of the shortest geodesic $\eta$: it has as many components as 
the order of $G$, since $G$ acts by freely permuting them. Choose one component 
of $\tilde\eta$, say $\eta_0$. Up to composing $\varphi$ with an element of $G$ 
we can assume that $\varphi(\eta_0)=\eta_0$. Note that $\eta_0$ is separated 
from the other components of $\tilde\eta$ by a unique lift $P_1\cup P_2$ of 
both pairs of pants bounding $H\setminus(\Lambda\cup\bigcup_{i=1}^3 C_i)$.

Let $X$ denote the lift of $H\setminus(\Lambda\cup\bigcup_{i=1}^3 C_i)$ bounded 
by $P_1\cup P_2$.

\medskip

\begin{Lemma}
We have $\varphi(P_1\cup P_2)\cap(P_1\cup P_2)\neq\emptyset$.
\end{Lemma}

\begin{proof}
Seeking a contradiction, assume that the intersection  is empty. Since $P_1\cup 
P_2$ separates and $\varphi(\eta_0)=\eta_0$, by volume reasons $P_1\cup P_2$ 
must separate $\varphi(P_1)$ from $\varphi(P_2)$. As every end of 
$\varphi(P_1)$ meets one end of $\varphi(P_2)$ along a cusp, $\varphi$  
stabilizes the set of three cusps $C_1\cup C_2\cup C_3$, joining $P_1$ and 
$P_2$. Thus $\varphi(P_1)$ (or $\varphi(P_2)$) is a totally geodesic embedded 
pant in $X$, with the  cusps $C_1\cup C_2\cup C_3$. By 
Lemma~\ref{Lemma:no_pants}, $\varphi(P_i)$ must be either $P_1$ or $P_2$, hence 
a contradiction.
\end{proof}

\begin{figure}[ht]
\begin{center}
\includegraphics[width=6cm]{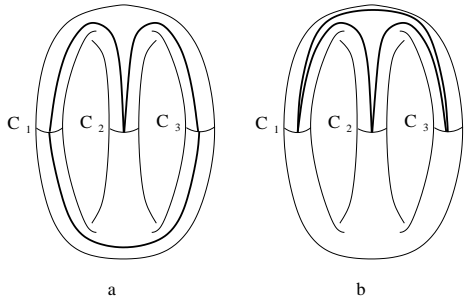}
\end{center}
\caption{The intersection $\varphi (P_j)\cap (P_1\cup P_2)$, viewed in
$\partial H=P_1\cup P_2 \cup \bigcup_{i=1}^3 C_i$.}
\label{triangleinhandle}
\end{figure}

\begin{Lemma}\label{claim:invariant boundary}
We have $\varphi(P_1\cup P_2)= P_1\cup P_2$.
\end{Lemma}
\begin{proof}
%
%
%
%
%
%
%
%
%
%
Since both $\varphi(P_1\cup P_2)$ and $P_1\cup P_2$ are totally geodesic, if 
they do not coincide, their intersection must be a union of disjoint simple 
geodesics, which are properly embedded. Figure~\ref{geod} shows the possible 
form of non self-intersecting proper geodesics, depending on whether the ends
are on the same or on different cusps. Observe that, since two distinct cusps 
cannot be sent by $\varphi$ to the same one, the geodesics in the intersection 
look the same in $P_1\cup P_2$ and in $\varphi(P_1\cup P_2)$.

\begin{figure}[ht]
\begin{center}
\includegraphics[width=9cm]{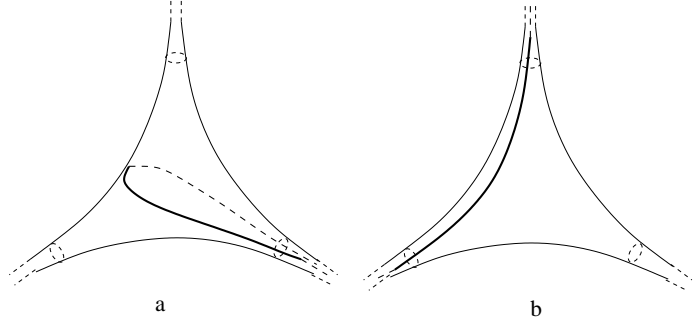}
\end{center}
\caption{The two types of geodesic intersections}
\label{geod}
\end{figure}

Assume that the intersection $P_i\cup \varphi(P_j)$ contains a geodesic as in 
Figure~\ref{geod}a, for some $i,j=1,2$. Notice that the intersection 
$\varphi(P_j)\cap (P_1\cup P_2)$ does not contain more geodesics, because 
$P_1\cup P_2$ separates. Consider now the half pant of $\varphi(P_j)$ contained 
inside $X$: it contains a cusp which ends either on a cusp inside $X$ or on a 
cusp corresponding to a $C_i$ different from that on which the geodesic of the 
intersection ends. In both cases we can find an annulus, properly embedded in 
$X$ and joining two cusps. The annulus is illustrated in Figure~\ref{inter}a: 
the shaded region is a disc contained in a section of the rank $1$ cusp. Since 
$X$ is hyperbolic, the two cusps which support the annulus, must in fact be the 
same and the annulus must be boundary parallel. Since totally geodesic parallel 
surfaces must be the same, it follows that $P_i=\varphi(P_j)$.


\begin{figure}[ht]
\begin{center}
\includegraphics[width=9cm]{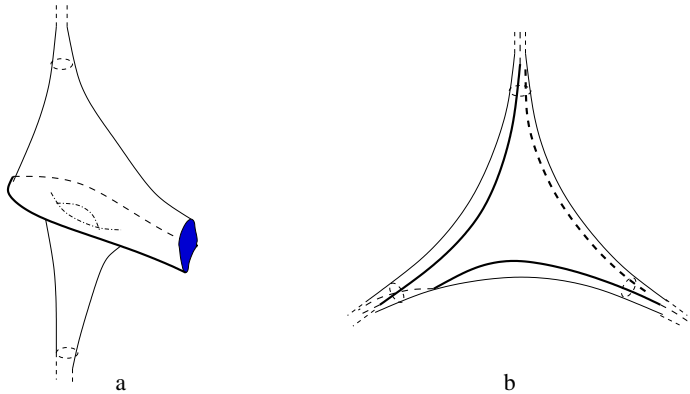}
\end{center}
\caption{An annulus and two ideal triangles}
\label{inter}
\end{figure}

We can thus assume that the intersection $\varphi(P_j)\cap (P_1\cup P_2)$
contains a geodesic of the type pictured in Figure~\ref{geod}b. Since 
$P_1\cup P_2$ separates, $\varphi(P_j)\cap (P_1\cup P_2)$ must contain two 
other geodesics of the same type that divide $\varphi(P_j)$ into two triangles, 
as illustrated in Figure~\ref{inter}b. If the three geodesics are not contained 
in a single pant $P_i$ (but in the union $P_1\cup P_2$), they give an essential 
curve in the genus two surface $\partial H=P_1\cup P_2\cup\bigcup_{i=1}^3 C_i$, 
which is the boundary of a compressing disc (half of $\varphi (P_j)$) (cf.\ 
Figure~\ref{triangleinhandle}a). Hence we may assume that the geodesics are 
contained in one of the pairs of pants $P_i$ and therefore they cut it into two 
ideal triangles, as in Figure~\ref{triangleinhandle}b. The union of half of 
$\varphi(P_j)$ and half of $P_i$, along the three geodesics and three segments 
in $C_1\cup C_2\cup C_3$, gives an embedded $2$-sphere. Irreducibility gives a 
parallelism between triangles in $P_i$ and $\varphi(P_j)$, hence 
$P_i=\varphi(P_j)$.

Again for volume reasons, we see that the remaining pairs of pants must 
intersect. Repeating the argument once more we reach the desired conclusion.
\end{proof}

\bigskip

Lemma~\ref{claim:invariant boundary} implies that $\varphi$ induces an isometry 
of $X=H\setminus(\Lambda\cup\bigcup_{i=1}^3 C_i)$. By 
Lemma~\ref{claim:asymmetric hyperbolic link} the restriction of $\varphi$ to 
$X$ is the identity. Since $X$ has nonempty interior in $\sph\setminus L$ and 
$\varphi$ is an isometry, we deduce that $\varphi$ itself is the identity and 
the desired contradiction follows. This finishes the proof of 
Theorem~\ref{thm}.
\qed

\section{Rigidifying actions on homology spheres}
\label{section:homology}

With the previous construction in mind we can now prove the corollary.

\medskip

\noindent \emph{Proof of Corollary~\ref{cor}}.
Consider a rational homology sphere $M$ on which $G$ acts freely \cite{CL}.
Repeating the previous construction, we can find a hyperbolic link $\L$ in $M$ 
which is $G$ invariant and such that $G$ is precisely the full 
orientation-preserving isometry group of the exterior of $\L$. The idea now is 
to do surgery in a $G$-equivariant way, so that we still have a rational 
homology sphere, and that the cores of the surgery tori are the shortest
curves, so that the $G$-orbits of those curves are invariant by any isometry.

To do this surgery, we must be able to choose a $G$-equivariant 
meridian-longitude system $(\mu,\lambda)$ on each peripheral torus (i.e.\ so 
that the image of $\lambda$ in $H^1(M\setminus \L;\QQ)$ is trivial). We specify 
how to adapt the construction of Theorem~\ref{thm}. First we remove a 
handlebody $H$ from the quotient $M/G$ and consider a hyperbolic anannular knot 
$\K$ in $M/G\setminus H$ in a trivial homotopy class, see \cite{My2}, so that it 
bounds a singular disc. This singular disc lifts to a family of $G$-equivariant 
singular discs, hence defining longitudes for the lifts of $\K$ in $M$.

The 
same construction with singular discs must be applied to the knots we remove 
from the interior of the handlebody $H$, but we need to justify that it is compatible 
with equivariance and the fact that we remove several curves:
\begin{itemize}
\item 
Recall from the proof of Lemma~\ref{claim:asymmetric hyperbolic link} that 
we remove a trivial knot $\gamma$ from the quotient
$\mathcal O= (H\setminus \bigcup_{i=1}^3C_i) \big/ \ZZ/2\times \mathbf D_3$, we lift
it   $\tilde \gamma\subset H\setminus\bigcup_{i=1}^3C_i$ and then we remove an equivariant hyperbolic 
anannular link $\Lambda_0$
from $H\setminus (\bigcup_{i=1}^3C_i\cup\tilde\gamma)$. We claim that  
$\Lambda_0$ can be choosen to project to a homotopically
trivial knot  in $\mathcal O\setminus \gamma$. To do that, we choose in $H$ a fundamental domain for the action
of $\ZZ/2\times \mathbf D_3$, and choose a homotopically trivial knot in $\mathcal O\setminus \gamma$ that 
crosses transversally each side of the fundamental domain at least twice (so that the punctured sides of the 
fundamental domain have negative Euler characteristic).
This gives a family of arcs in the fundamental domain, that can be homotoped relative to the boundary to   
a submanifold with hyperbolic and ananular exterior, by the main theorem in \cite{My2}. By the gluing
lemma of Myers, \cite[Lemma~2.1]{My2}, the pieces of the fundamental domain match to give the link
$\Lambda_0$ with the required properties.
\item 
We justify the compatibility of singular discs when we remove more than one curve, that 
can intersect such a singular disc, or the meridian discs bounding the curves 
$C_i$ (that we removed from $\partial H$). In this case we tube the discs along 
the curves we remove, in order to get a disjoint surface. This tubing is 
possible because, at each step, all curves we remove are homotopically trivial. 
\end{itemize}
Finally, notice that the perturbation argument of Proposition~\ref{prop}, which is not 
used here, would be a problem for the existence of longitudes.


We now perform hyperbolic Dehn surgery with filling meridians $n\lambda+\mu$  
on the components of $\L$, $n \gg 1$, in such a way that the following 
requirements are satisfied:

\begin{itemize}

\item the surgery is $G$-equivariant;

\item the geodesics obtained after surgery have pairwise distinct lengths if 
they belong to different $G$-orbits;

\item the lengths of these geodesics are the shortest ones.


\end{itemize}

Notice that the slopes are chosen so that the resulting manifold is a rational
homology sphere. Since we choose  homotopically trivial curves in the quotient, 
the $G$ action permutes the componets of the link, hence $G$ acts 
freely on the surgered manifold. 
The geodesics are then 
chosen to be of shortest lengths to ensure that all isometries must preserve 
the image of $\L$ after surgery, so that they must induce isometries of the 
complement of $\L$. The conclusion now follows at once.
\qed

\paragraph*{Acknowledgments} 
 The first   author wishes to 
 thank
the Departament de Matem\`atiques of the Universitat Aut\`onoma for hospitality
while this work was carried out.

\bigskip

\begin{footnotesize}
\textsc{IMB - UMR 5584,
	Universit\'e de Bourgogne,
	BP 47870,
	9 av. Alain Savary}
	
\textsc{21000 Dijon CEDEX
	France}
	
{paoluzzi@u-bourgogne.fr}

\medskip

\textsc{Departament de Matem\`atiques, Universitat Aut\`onoma de Barcelona.}

\textsc{08193 Bellaterra, Spain}

{porti@mat.uab.es}

\end{footnotesize}

\end{document}